# Where do statistical models come from? Revisiting the problem of specification

Aris Spanos*[1]

*Virginia Polytechnic Institute and State University*

**Abstract:** R. A. Fisher founded modern statistical inference in 1922 and identified its fundamental problems to be: *specification, estimation* and *distribution*. Since then the problem of *statistical model specification* has received scant attention in the statistics literature. The paper traces the history of statistical model specification, focusing primarily on pioneers like Fisher, Neyman, and more recently Lehmann and Cox, and attempts a synthesis of their views in the context of the Probabilistic Reduction (PR) approach. As argued by Lehmann [11], a major stumbling block for a general approach to statistical model specification has been the delineation of the appropriate role for substantive subject matter information. The PR approach demarcates the interrelated but complemenatry roles of *substantive* and *statistical information* summarized *ab initio* in the form of a *structural* and a *statistical model*, respectively. In an attempt to preserve the integrity of both sources of information, as well as to ensure the reliability of their fusing, a purely probabilistic construal of statistical models is advocated. This probabilistic construal is then used to shed light on a number of issues relating to specification, including the role of preliminary data analysis, structural vs. statistical models, model specification vs. model selection, statistical vs. substantive adequacy and model validation.

## 1. Introduction

The current approach to statistics, interpreted broadly as 'probability-based data modeling and inference', has its roots going back to the early 19th century, but it was given its current formulation by R. A. Fisher [5]. He identified the fundamental problems of statistics to be: *specification, estimation* and *distribution*. Despite its importance, the question of *specification,* 'where do statistical models come from?' received only scant attention in the statistics literature; see Lehmann [11].

The cornerstone of modern statistics is the notion of a *statistical model* whose meaning and role have changed and evolved along with that of statistical modeling itself over the last two centuries. Adopting a retrospective view, a statistical model is defined to be an internally consistent set of probabilistic assumptions aiming to provide an 'idealized' probabilistic description of the stochastic mechanism that gave rise to the observed data $\mathbf{x} := (x_1, x_2, \ldots, x_n)$. The quintessential statistical model is the *simple Normal model*, comprising a statistical Generating Mechanism (GM):

$$(1.1) \qquad X_k = \mu + u_k, \ k \in \mathbb{N} := \{1, 2, \ldots n, \ldots\}$$

*I'm most grateful to Erich Lehmann, Deborah G. Mayo, Javier Rojo and an anonymous referee for valuable suggestions and comments on an earlier draft of the paper.
[1]Department of Economics, Virginia Polytechnic Institute, and State University, Blacksburg, VA 24061, e-mail: aris@vt.edu

together with the probabilistic assumptions:

(1.2) $$X_k \sim \mathsf{NIID}(\mu, \sigma^2), \ k \in \mathbb{N},$$

where $X_k \sim$ NIID stands for *Normal, Independent and Identically Distributed*. The nature of a statistical model will be discussed in section 3, but as a prelude to that, it is important to emphasize that it is specified exclusively in terms of probabilistic concepts that can be related directly to the *joint distribution* of the *observable* stochastic process $\{X_k, \ k \in \mathbb{N}\}$. This is in contrast to other forms of models that play a role in statistics, such as *structural* (explanatory, substantive), which are based on substantive subject matter information and are specified in terms of theory concepts.

The motivation for such a purely probabilistic construal of statistical models arises from an attempt to circumvent some of the difficulties for a general approach to statistical modeling. These difficulties were raised by early pioneers like Fisher [5]–[7] and Neyman [17]–[26], and discussed extensively by Lehmann [11] and Cox [1]. The main difficulty, as articulated by Lehmann [11], concerns the role of *substantive subject matter* information. His discussion suggests that if statistical model specification requires such information at the outset, then any attempt to provide a general approach to statistical modeling is unattainable. His main conclusion is that, despite the untenability of a general approach, statistical theory has a contribution to make in model specification by extending and improving: (a) the reservoir of models, (b) the model selection procedures, and (c) the different types of models.

In this paper it is argued that Lehmann's case concerning (a)–(c) can be strengthened and extended by adopting a purely probabilistic construal of statistical models and placing statistical modeling in a broader framework which allows for fusing statistical and substantive information in a way which does not compromise the integrity of either. Substantive subject matter information emanating from the theory, and statistical information reflecting the probabilistic structure of the data, need to be viewed as bringing to the table different but complementary information. The *Probabilistic Reduction* (PR) approach offers such a modeling framework by integrating several innovations in Neyman's writings into Fisher's initial framework with a view to address a number of modeling problems, including the role of preliminary data analysis, structural vs. statistical models, model specification vs. model selection, statistical vs. substantive adequacy and model validation. Due to space limitations the picture painted in this paper will be dominated by broad brush strokes with very few details; see Spanos [31]–[42] for further discussion.

### *1.1. Substantive vs. statistical information*

Empirical modeling in the social and physical sciences involves an intricate blending of *substantive subject matter* and *statistical information*. Many aspects of empirical modeling implicate both sources of information in a variety of functions, and others involve one or the other, more or less separately. For instance, the development of *structural* (explanatory) *models* is primarily based on *substantive information* and it is concerned with the mathematization of theories to give rise to theory models, which are amenable to empirical analysis; that activity, by its very nature, cannot be separated from the disciplines in question. On the other hand, certain aspects of empirical modeling, which focus on *statistical information* and are concerned with the nature and use of statistical models, can form a body of knowledge which is shared by all fields that use data in their modeling. This is the body of knowledge that *statistics* can claim as its subject matter and develop it with only one eye on new problems/issues raised by empirical modeling in other disciplines. This ensures that statistics is not subordinated to the other applied fields, but remains a separate discipline which provides, maintains and extends/develops the common foundation and overarching framework for empirical modeling.

To be more specific, *statistical model specification* refers to the choice of a model (parameterization) arising from the probabilistic structure of a stochastic process $\{X_k,\ k\in\mathbb{N}\}$ that would render the data in question $\mathbf{x}:=(x_1, x_2, \ldots, x_n)$ a *truly typical realization* thereof. This perspective on the data is referred to as the Fisher–Neyman probabilistic perspective for reasons that will become apparent in section 2. When one specifies the simple Normal model (1.1), the only thing that matters from the statistical specification perspective is whether the data $\mathbf{x}$ can be realistically viewed as a truly typical realization of the process $\{X_k,\ k\in\mathbb{N}\}$ assumed to be NIID, devoid of any substantive information. A model is said to be *statistical adequate* when the assumptions constituting the statistical model in question, NIID, are valid for the data $\mathbf{x}$ in question. Statistical adequacy can be assessed qualitatively using analogical reasoning in conjunction with data graphs (t-plots, P-P plots etc.), as well as quantitatively by testing the assumptions constituting the statistical model using probative Mis-Specification (M-S) tests; see Spanos [36].

It is argued that certain aspects of statistical modeling, such as statistical model specification, the use of graphical techniques, M-S testing and respecification, together with optimal inference procedures (estimation, testing and prediction), can be developed *generically* by viewing data $\mathbf{x}$ as a realization of a (nameless) stochastic process $\{X_k,\ t\in\mathbb{N}\}$. All these aspects of empirical modeling revolve around a central axis we call a *statistical model*. Such models can be viewed as *canonical models,* in the sense used by Mayo [12], which are developed without any reference to substantive subject matter information, and can be used equally in physics, biology, economics and psychology. Such canonical models and the associated statistical modeling and inference belong to the *realm of statistics*. Such a view will broaden the scope of modern statistics by integrating preliminary data analysis, statistical model specification, M-S testing and respecification into the current textbook discourses; see Cox and Hinkley [2], Lehmann [10].

On the other hand the question of *substantive adequacy*, i.e. whether a structural model adequately captures the main features of the *actual Data Generating Mechanism* (DGM) giving rise to data $\mathbf{x}$, cannot be addressed in a generic way because it concerns the bridge between the particular model and the phenomenon of interest. Even in this case, however, assessing substantive adequacy will take the form of applying statistical procedures within an embedding statistical model. Moreover, for the error probing to be *reliable* one needs to ensure that the embedding model is statistically adequate; it captures all the statistical systematic information (Spanos, [41]). In this sense, substantive subject matter information (which might range from vary vague to highly informative) constitutes important supplementary information which, under *statistical* and *substantive adequacy*, enhances the explanatory and predictive power of statistical models.

In the spirit of Lehmann [11], models in this paper are classified into:

(a) *statistical* (empirical, descriptive, interpolatory formulae, data models), and
(b) *structural* (explanatory, theoretical, mechanistic, substantive).

The harmonious integration of these two sources of information gives rise to an

(c) *empirical model*; the term is not equivalent to that in Lehmann [11].

In Section 2, the paper traces the development of ideas, issues and problems surrounding statistical model specification from Karl Pearson [27] to Lehmann [11], with particular emphasis on the perspectives of Fisher and Neyman. Some of the ideas and modeling suggestions of these pioneers are synthesized in Section 3 in the form of the PR modeling framework. Kepler's first law of planetary motion is used to illustrate some of the concepts and ideas. The PR perspective is then used to shed light on certain issues raised by Lehmann [11] and Cox [1].

## 2. 20th century statistics

### 2.1. Early debates: description vs. induction

Before Fisher, the notion of a statistical model was both vague and implicit in data modeling, with its role primarily confined to the *description* of the distributional properties of *the data in hand* using the histogram and the first few sample moments. A crucial problem with the application of *descriptive statistics* in the late 19th century was that statisticians would often claim generality beyond the data in hand for their inferences. This is well-articulated by Mills [16]:

"In approaching this subject [statistics] we must first make clear the distinction between *statistical description* and *statistical induction.* By employing the methods of statistics it is possible, as we have seen, to describe succinctly a mass of quantitative data." … "In so far as the results are confined to the cases actually studied, these various statistical measures are merely devices for describing certain features of a distribution, or certain relationships. Within these limits the measures may be used to perfect confidence, as accurate descriptions of the given characteristics. But when we seek to extend these results, to generalize the conclusions, to apply them to cases not included in the original study, a quite new set of problems is faced." (p. 548-9)

Mills [16] went on to discuss the 'inherent assumptions' necessary for the validity of *statistical induction:*

"… in the larger population to which this result is to be applied, there exists a uniformity with respect to the characteristic or relation we have measured" …, and "… the sample from which our first results were derived is thoroughly representative of the entire population to which the results are to be applied." (pp. 550-2).

The fine line between statistical description and statistical induction was nebulous until the 1920s, for several reasons. *First*, "No distinction was drawn between a sample and the population, and what was calculated from the sample was attributed to the population." (Rao, [29], p. 35). *Second*, it was thought that the inherent assumptions for the validity of statistical induction are *not empirically verifiable*; see Mills [16], p. 551). *Third*, there was a widespread belief, exemplified in the first quotation from Mills, that statistical description does not require any assumptions. It is well-known today that there is no such thing as a meaningful summary of the data that does not involve any implicit assumptions; see Neyman [21]. For instance, the arithmetic average of a trending time series represents no meaningful feature of the underlying 'population'.

## 2.2. Karl Pearson

Karl Pearson was able to take descriptive statistics to a higher level of sophistication by proposing the 'graduation (smoothing) of histograms' into 'frequency curves'; see Pearson [27]. This, however, introduced additional fuzziness into the distinction between statistical description vs. induction because the frequency curves were the precursors to the density functions; one of the crucial components of a statistical model introduced by Fisher [5] providing the foundation of statistical induction. The statistical modeling procedure advocated by Pearson, however, was very different from that introduced by Fisher.

For Karl Pearson statistical modeling would begin with data $\mathbf{x} := (x_1, x_2, \ldots, x_n)$ in search of a descriptive model which would be in the form of a frequency curve $f(x)$, chosen from the Pearson family $f(x; \theta)$, $\theta := (a, b_0, b_1, b_2)$, after applying the method of moments to obtain $\widehat{\theta}$ (see Pearson, [27]). Viewed from today's perspective, the solution $\widehat{\theta}$, would deal with two different statistical problems simultaneously, (a) specification (the choice a descriptive model $f(x; \widehat{\theta})$) and (b) *estimation* of $\theta$ using $\widehat{\theta}$. $f(x; \widehat{\theta})$ can subsequently be used to draw inferences beyond the original data $\mathbf{x}$.

Pearson's view of statistical induction, as late as 1920, was that of *induction by enumeration* which relies on both *prior* distributions and the stability of relative frequencies; see Pearson [28], p. 1.

## 2.3. R. A. Fisher

One of Fisher's most remarkable but least appreciated achievements, was to initiate the recasting of the form of *statistical induction* into its modern variant. Instead of starting with data $\mathbf{x}$ in search of a descriptive model, he would interpret the data as a truly *representative sample* from a pre-specified 'hypothetical infinite population'. This might seem like a trivial re-arrangement of Pearson's procedure, but in fact it constitutes a complete recasting of the problem of statistical induction, with the notion of a *parameteric statistical model* delimiting its premises.

Fisher's first clear statement of this major change from the then prevailing modeling process is given in his classic 1922 paper:

"... the object of statistical methods is the reduction of data. A quantity of data, which usually by its mere bulk is incapable of entering the mind, is to be replaced by relatively few quantities which shall adequately represent the whole, or which, in other words, shall contain as much as possible, ideally the whole, of the relevant information contained in the original data. This object is accomplished by constructing a hypothetical infinite population, of which the actual data are regarded as constituting a sample. The law of distribution of this hypothetical population is specified by relatively few parameters, which are sufficient to describe it exhaustively in respect of all qualities under discussion." ([5], p. 311)

Fisher goes on to elaborate on the modeling process itself: "The problems which arise in reduction of data may be conveniently divided into three types: (1) **Problems of Specification**. These arise in the choice of the mathematical form of the population. (2) **Problems of Estimation**. (3) **Problems of Distribution**.

It will be clear that when we know (1) what parameters are required to specify the population from which the sample is drawn, (2) how best to calculate from the sample estimates of these parameters, and (3) the exact form of the distribution, in different samples, of our derived statistics, then the theoretical aspect of the treatment of any particular body of data has been completely elucidated." (p. 313-4)

One can summarize Fisher's view of the statistical modeling process as follows. The process begins with a *prespecified parametric statistical model* $\mathcal{M}$ ('a hypothetical infinite population'), chosen so as to ensure that the observed data **x** are viewed as a *truly representative sample* from that 'population':

"The postulate of randomness thus resolves itself into the question, "Of what population is this a random sample?" which must frequently be asked by every practical statistician." ([5], p. 313)

Fisher was fully aware of the fact that the specification of a statistical model premises all forms of statistical inference. Once $\mathcal{M}$ was specified, the original uncertainty relating to the 'population' was reduced to uncertainty concerning the unknown parameter(s) $\theta$, associated with $\mathcal{M}$. In Fisher's set up, the parameter(s) $\theta$, are *unknown constants* and become the focus of inference. The problems of 'estimation' and 'distribution' revolve around $\theta$.

Fisher went on to elaborate further on the 'problems of specification': "As regards problems of specification, these are entirely a matter for the practical statistician, for those cases where the qualitative nature of the hypothetical population is known do not involve any problems of this type. In other cases we may know by experience what forms are likely to be suitable, and *the adequacy of our choice may be tested a posteriori*. We must confine ourselves to those forms which we know how to handle, or for which any tables which may be necessary have been constructed. More or less elaborate form will be suitable according to the volume of the data." (p. 314) [emphasis added]

Based primarily on the above quoted passage, Lehmann's [11] assessment of Fisher's view on specification is summarized as follows: "Fisher's statement implies that in his view there can be no theory of modeling, no general modeling strategies, but that instead each problem must be considered entirely on its own merits. He does not appear to have revised his opinion later... Actually, following this uncompromisingly negative statement, Fisher unbends slightly and offers two general suggestions concerning model building: (a) "We must confine ourselves to those forms which we know how to handle," and (b) "More or less elaborate forms will be suitable according to the volume of the data."" (p. 160-1).

Lehmann's interpretation is clearly warranted, but Fisher's view of specification has some additional dimensions that need to be brought out. The original choice of a statistical model may be guided by simplicity and experience, but as Fisher emphasizes "the adequacy of our choice may be tested *a posteriori*." What comes after the above quotation is particularly interesting to be quoted in full: "Evidently these are considerations the nature of which may change greatly during the work of a single generation. We may instance the development by Pearson of a very extensive system of skew curves, the elaboration of a method of calculating their parameters, and the preparation of the necessary tables, a body of work which has enormously extended the power of modern statistical practice, and which has been, by pertinacity and inspiration alike, practically the work of a single man. Nor is the introduction of the Pearsonian system of frequency curves the only contribution which their author has made to the solution of problems of specification: of even greater importance is the introduction of an objective criterion of goodness of fit. For empirical as the specification of the hypothetical population may be, *this empiricism* *is cleared of its dangers if we can apply a rigorous and objective test of the adequacy with which the proposed population represents the whole of the available facts.* Once a statistic suitable for applying such a test, has been chosen, the exact form of its distribution in random samples must be investigated, in order that we may evaluate the probability that a worse fit should be obtained from a random sample of a population of the type considered. The possibility of developing complete and self-contained tests of goodness of fit deserves very careful consideration, *since therein lies our justification for the free use which is made of empirical frequency formulae.* Problems of distribution of great mathematical difficulty have to be faced in this direction." (p. 314) [emphasis (in italic) added]

In this quotation Fisher emphasizes the empirical dimension of the specification problem, and elaborates on testing the assumptions of the model, lavishing Karl Pearson with more praise for developing the *goodness of fit test* than for his family of densities. He clearly views this test as a primary tool for assessing the validity of the original specification (misspecification testing). He even warns the reader of the potentially complicated sampling theory required for such form of testing. Indeed, most of the tests he discusses in chapters 3 and 4 of his 1925 book [6] are *misspecification tests*: tests of departures from Normality, Independence and Homogeneity. Fisher emphasizes the fact that the reliability of every form of inference depend crucially on the validity of the statistical model postulated. The premises of statistical induction in Fisher's sense no longer rely on prior assumptions of 'ignorance', but on testable probabilistic assumptions which concern the observed data; this was a major departure from Pearson's form of enumerative induction relying on prior distributions.

A more complete version of the three problems of the 'reduction of data' is repeated in Fisher's 1925 book [6], which is worth quoting in full with the major additions indicated in italic: "The problems which arise in the reduction of data may thus conveniently be divided into three types:

(i) Problems of Specification, which arise in the choice of the mathematical form of the population. *This is not arbitrary, but requires an understanding of the way in which the data are supposed to, or did in fact, originate. Its further discussion depends on such fields as the theory of Sample Survey, or that of Experimental Design.*

(ii) *When the specification has been obtained, problems of Estimation arise.* These involve the choice among the methods of calculating, from our sample, statistics fit to estimate the unknown parameters of the population.

(iii) Problems of *Distribution* include the mathematical deduction of the exact nature of the distributions in random samples of our estimates of the parameters, *and of the other statistics designed to test the validity of our specification (tests of Goodness of Fit).*" (see ibid. p. 8)

In (i) Fisher makes a clear reference to the *actual Data Generating Mechanism* (DGM), which might often involve specialized knowledge beyond statistics. His view of specification, however, is narrowed down by his focus on data from 'sample surveys' and 'experimental design', where the gap between the actual DGM and the statistical model is not sizeable. This might explain his claim that: "... for those cases where the qualitative nature of the hypothetical population is known do not involve any problems of this type." In his 1935 book, Fisher states that: "Statistical procedure and experimental design are only two aspects of the same whole, and that whole comprises all the logical requirements of the complete process of adding to natural knowledge by experimentation" (p. 3)

In (iii) Fisher adds the derivation of the sampling distributions of misspecification tests as part of the 'problems of distribution'.

In summary, Fisher's view of specification, as a facet of modeling providing the foundation and the overarching framework for *statistical induction*, was a radical departure from Karl Pearson's view of the problem. By interpreting the observed data as 'truly representative' of a prespecified statistical model, Fisher initiated the recasting of *statistical induction* and rendered its premises testable. By ascertaining *statistical adequacy*, using *misspecification tests*, the modeler can ensure the reliability of inductive inference. In addition, his pivotal contributions to the 'problems of Estimation and Distribution', in the form of finite sampling distributions for estimators and test statistics, shifted the emphasis in statistical induction, from enumerative induction and its reliance on asymptotic arguments, to 'reliable procedures' based on finite sample 'ascertainable error probabilities':

"In order to assert that a natural phenomenon is experimentally demonstrable we need, not an isolated record, but a reliable method of procedure. In relation to the test of significance, we may say that a phenomenon is experimentally demonstrable when we know how to conduct an experiment which will rarely fail to give us a statistically significant result." (Fisher [7], p. 14)

This constitutes a clear description of *inductive inference* based on ascertainable error probabilities, under the 'control' of the experimenter, used to assess the 'optimality' of inference procedures. Fisher was the first to realize that for precise (finite sample) 'error probabilities', to be used for calibrating statistical induction, one needs a complete model specification including a distribution assumption. Fisher's most enduring contribution is his devising a general way to 'operationalize' the errors for statistical induction by embedding the *material experiment* into a *statistical model* and define the frequentist error probabilities in the context of the latter. These statistical error probabilities provide a measure of the 'trustworthiness' of the inference procedure: how often it will give rise to true inferences concerning the underlying DGM. That is, the inference is reached by an inductive procedure which, with high probability, will reach true conclusions from true (or approximately true) premises (statistical model). This is in contrast to *induction by enumeration* where the focus is on observed 'events' and not on the 'process' generating the data.

In relation to this, C. S. Peirce put forward a similar view of quantitative induction, almost half a century earlier. This view of statistical induction, was called the error statistical approach by Mayo [12], who has formalized and extended it to include a post-data evaluation of inference in the form of severe testing. Severe testing can be used to address chronic problems associated with Neyman-Pearson testing, including the classic fallacies of acceptance and rejection; see Mayo and Spanos [14].

### *2.4. Neyman*

According to Lehmann [11], Neyman's views on the theory of statistical modeling had three distinct features:

"1. Models of complex phenomena are constructed by combining simple building blocks which, "partly through experience and partly through imagination, appear to us familiar, and therefore, simple." ...

2. An important contribution to the theory of modeling is Neyman's distinction between two types of models: "interpolatory formulae" on the one hand and "explanatory models" on the other. The latter try to provide an explanation of the mechanism underlying the observed phenomena; Mendelian inheritance was Neyman's favorite example. On the other hand an interpolatory formula is based on a convenient and flexible family of distributions or models given a priori, for example the Pearson curves, one of which is selected as providing the best fit to the data. ...

3. The last comment of Neyman's we mention here is that to develop a "genuine explanatory theory" requires substantial knowledge of the scientific background of the problem." (p. 161)

Lehmann's first hand knowledge of Neyman's views on modeling is particularly enlightening. It is clear that Neyman adopted, adapted and extended Fisher's view of statistical modeling. What is especially important for our purposes is to bring out both the similarities as well as the subtle differences with Fisher's view.

Neyman and Pearson [26] built their hypothesis testing procedure in the context of Fisher's approach to statistical modeling and inference, with the notion of a prespecified parametric statistical model providing the cornerstone of the whole inferential edifice. Due primarily to Neyman's experience with empirical modeling in a number of applied fields, including genetics, agriculture, epidemiology and astronomy, his view of statistical models, evolved beyond Fisher's 'infinite populations' in the 1930s into *frequentist* 'chance mechanisms' in the 1950s:

"(ii) Guessing and then verifying the 'chance mechanism', the repeated operations of which produces the observed frequencies. This is a problem of 'frequentist probability theory'. Occasionally, this step is labelled 'model building'. Naturally, the guessed chance mechanism is hypothetical." (Neyman [25], p. 99)

In this quotation we can see a clear statement concerning the nature of specification. Neyman [18] describes statistical modeling as follows: "The application of the theory involves the following steps:

(i) If we wish to treat certain phenomena by means of the theory of probability we must find some element of these phenomena that could be considered as random, following the law of large numbers. This involves a construction of a mathematical model of the phenomena involving one or more probability sets.

(ii) The mathematical model is found satisfactory, or not. This must be checked by observation.

(iii) If the mathematical model is found satisfactory, then it may be used for deductions concerning phenomena to be observed in the future." (ibid., p. 27)

In this quotation Neyman in (i) demarcates the domain of statistical modeling to *stochastic phenomena*: observed phenomena which exhibit *chance regularity patterns*, and considers statistical (mathematical) models as probabilistic constructs. He also emphasizes the reliance of frequentist inductive inference on the long-run stability of relative frequencies. Like Fisher, he emphasizes in (ii) the testing of the assumptions comprising the statistical model in order to ensure its adequacy. In (iii) he clearly indicates that statistical adequacy is a necessary condition for any inductive inference. This is because the 'error probabilities', in terms of which the optimality of inference is defined, depend crucially on the validity of the model:

"... any statement regarding the performance of a statistical test depends upon the postulate that the observable random variables are random variables and posses the properties specified in the definition of the set $\Omega$ of the admissible simple hypotheses." (Neyman [17], p. 289)

A crucial implication of this is that when the statistical model is misspecified, the actual error probabilities, in terms of which 'optimal' inference procedures are chosen, are likely to be very different from the nominal ones, leading to unreliable inferences; see Spanos [40].

Neyman's experience with modeling *observational data* led him to take statistical modeling a step further and consider the question of *respecifying* the original model whenever it turns out to be inappropriate (statistically inadequate): "Broadly, the methods of bringing about an agreement between the predictions of statistical theory and observations may be classified under two headings:(a) Adaptation of the statistical theory to the enforced circumstances of observation. (b) Adaptation of the experimental technique to the postulates of the theory. The situations referred to in (a) are those in which the observable random variables are largely outside the control of the experimenter or observer." ([17], p. 291)

Neyman goes on to give an example of (a) from his own applied research on the effectiveness of insecticides where the Poisson model was found to be inappropriate: "Therefore, if the statistical tests based on the hypothesis that the variables follow the Poisson Law are not applicable, the only way out of the difficulty is to modify or adapt the theory to the enforced circumstances of experimentation." (ibid., p. 292)

In relation to (b) Neyman continues (ibid., p. 292): "In many cases, particularly in laboratory experimentation, the nature of the observable random variables is much under the control of the experimenter, and here it is usual to adapt the experimental techniques so that it agrees with the assumptions of the theory."

He goes on to give due credit to Fisher for introducing the crucially important technique of *randomization* and discuss its application to the 'lady tasting tea' experiment. Arguably, Neyman's most important extension of Fisher's specification facet of statistical modeling, was his underscoring of *the gap between a statistical model and the phenomena of interest*:

"...it is my strong opinion that no mathematical theory refers exactly to happenings in the outside world and that any application requires a solid bridge over an abyss. The construction of such a bridge consists first, in explaining in what sense the mathematical model provided by the theory is expected to "correspond" to certain actual happenings and second, in checking empirically whether or not the correspondence is satisfactory." ([18], p. 42)

He emphasizes the bridging of the gap between a statistical model and the observable phenomenon of interest, arguing that, beyond *statistical adequacy*, one needs to ensure *substantive adequacy*: the accord between the statistical model and 'reality' must also be adequate: "Since in many instances, the phenomena rather than their models are the subject of scientific interest, the transfer to the phenomena of an inductive inference reached within the model must be something like this: granting that the model $M$ of phenomena $P$ is adequate (or valid, of satisfactory, etc.) the conclusion reached within $M$ applies to $P$." (Neyman [19], p. 17)

In a purposeful attempt to bridge this gap, Neyman distinguished between a *statistical model* (interpolatory formula) and a *structural model* (see especially Neyman [24], p. 3360), and raised the important issue of *identification* in Neyman [23]: "This particular finding by Polya demonstrated a phenomenon which was unanticipated – two radically different stochastic mechanisms can produce identical distributions of the same variable $X$! Thus, the study of this distribution cannot answer the question which of the two mechanisms is actually operating. " ([23], p. 158)

In summary, Neyman's views on statistical modeling elucidated and extended that of Fisher's in several important respects: (a) Viewing statistical models primarily as 'chance mechanisms'. (b) Articulating fully the role of 'error probabilities' in assessing the optimality of inference methods. (c) Elaborating on the issue of *re-specification* in the case of statistically *inadequate* models. (d) Emphasizing the gap between a *statistical model* and the *phenomenon of interest*. (e) Distinguishing between *structural* and *statistical models*. (f) Recognizing the problem of *Identification*.

### 2.5. Lehmann

Lehmann [11] considers the question of 'what contribution statistical theory can potentially make to model specification and construction'. He summarizes the views of both Fisher and Neyman on model specification and discusses the meagre subsequent literature on this issue. His primary conclusion is rather pessimistic: apart from some vague guiding principles, such as simplicity, imagination and the use of past experience, no general theory of modeling seems attainable: "This requirement [to develop a "genuine explanatory theory" requires substantial knowledge of the scientific background of the problem] is agreed on by all serious statisticians but it constitutes of course an obstacle to any general theory of modeling, and is likely a principal reason for Fisher's negative feeling concerning the possibility of such a theory." (Lehmann [11], p. 161)

Hence, Lehmann's source of pessimism stems from the fact that 'explanatory' models place a major component of model specification beyond the subject matter of the statistician: "An explanatory model, as is clear from the very nature of such models, requires detailed knowledge and understanding of the substantive situation that the model is to represent. On the other hand, an empirical model may be obtained from a family of models selected largely for convenience, on the basis solely of the data without much input from the underlying situation." (p. 164)

In his attempt to demarcate the potential role of statistics in a general theory of modeling, Lehmann [11], p. 163, discusses the difference in the basic objectives of the two types of models, arguing that: "Empirical models are used as a guide to action, often based on forecasts ... In contrast, explanatory models embody the search for the basic mechanism underlying the process being studied; they constitute an effort to achieve *understanding*."

In view of these, he goes on to pose a crucial question (Lehmann [11], p. 161-2): "Is applied statistics, and more particularly model building, an art, with each new case having to be treated from scratch, ..., completely on its own merits, or does theory have a contribution to make to this process?"

Lehmann suggests that one (indirect) way a statistician can contribute to the theory of modeling is via: "... the existence of a reservoir of models which are well understood and whose properties we know. Probability theory and statistics have provided us with a rich collection of such models." (p. 161)

Assuming the existence of a sizeable reservoir of models, the problem still remains 'how does one make a choice among these models?' Lehmann's view is that the current methods on model selection *do not* address this question:

"Procedures for choosing a model not from the vast storehouse mentioned in (2.1 Reservoir of Models) but from a much more narrowly defined class of models are discussed in the theory of model selection. A typical example is the choice of a regression model, for example of the best dimension in a nested class of such models. ... However, this view of model selection ignores a preliminary step: the specification of the class of models from which the selection is to be made." (p. 162)

This is a most insightful comment because a closer look at model selection procedures suggests that the problem of model specification is largely assumed away by commencing the procedure by assuming that the prespecified family of models includes the true model; see Spanos [42].

In addition to differences in their nature and basic objectives, Lehmann [11] argues that explanatory and empirical models pose very different problems for model validation: "The difference in the aims and nature of the two types of models [empirical and explanatory] implies very different attitudes toward checking their validity. Techniques such as goodness of fit test or cross validation serve the needs of checking an empirical model by determining whether the model provides an adequate fit for the data. Many different models could pass such a test, which reflects the fact that there is not a unique correct empirical model. On the other hand, ideally there is only one model which at the given level of abstraction and generality describes the mechanism or process in question. To check its accuracy requires identification of the details of the model and their functions and interrelations with the corresponding details of the real situation." (ibid. pp. 164-5)

Lehmann [11] concludes the paper on a more optimistic note by observing that statistical theory has an important role to play in model specification by extending and enhancing: (1) the reservoir of models, (2) the model selection procedures, as well as (3) utilizing different classifications of models. In particular, in addition to the subject matter, every model also has a 'chance regularity' dimension and probability theory can play a crucial role in 'capturing' this. This echoes Neyman [21], who recognized the problem posed by explanatory (stochastic) models, but suggested that probability theory does have a crucial role to play: "The problem of stochastic models is of prime interest but is taken over partly by the relevant substantive disciplines, such as astronomy, physics, biology, economics, etc., and partly by the theory of probability. In fact, the primary subject of the modern theory of probability may be described as the study of properties of particular chance mechanisms." (p. 447)

Lehmann's discussion of model specification suggests that *the* major stumbling block in the development of a general modeling procedure is the substantive knowledge, beyond the scope of statistics, called for by explanatory models; see also Cox and Wermuth [3]. To be fair, both Fisher and Neyman in their writings seemed to suggest that statistical model specification is based on an amalgam of substantive and statistical information.

Lehmann [11] provides a key to circumventing this stumbling block: "Examination of some of the classical examples of revolutionary science shows that the eventual explanatory model is often reached in stages, and that in the earlier efforts one may find models that are descriptive rather than fully explanatory. ... This is, for example, true of Kepler whose descriptive model (laws) of planetary motion precede Newton's explanatory model." (p. 166).

In this quotation, Lehmann acknowledges that a descriptive (statistical) model can have 'a life of its own', separate from substantive subject matter information. However, the question that arises is: 'what is such model a description of?' As argued in the next section, in the context of the Probabilistic Reduction (PR) framework, such a model provides a description of the systematic *statistical infor-* *mation* exhibited by data $\mathbf{Z} := (\mathbf{z}_1, \mathbf{z}_2, \ldots, \mathbf{z}_n)$. This raises another question 'how does the substantive information, when available, enter statistical modeling?' Usually substantive information enters the empirical modeling as restrictions on a statistical model, when the structural model, carrying the substantive information, is embedded into a statistical model. As argued next, when these restrictions are data-acceptable, assessed in the context of a statistically adequate model, they give rise to an *empirical model* (see Spanos, [31]), which is both statistically as well as substantively meaningful.

## 3. The Probabilistic Reduction (PR) Approach

The foundations and overarching framework of the *PR approach* (Spanos, [31]–[42]) has been greatly influenced by Fisher's recasting of statistical induction based on the notion of a statistical model, and calibrated in terms of frequentist error probabilities, Neyman's extensions of Fisher's paradigm to the modeling of observational data, and Kolmogorov's crucial contributions to the theory of stochastic processes. The emphasis is placed on learning from data about observable phenomena, and on actively encouraging thorough probing of the different ways an inference might be in error, by localizing the error probing in the context of different models; see Mayo [12]. Although the broader problem of bridging the gap between theory and data using a sequence of interrelated models (see Spanos, [31], p. 21) is beyond the scope of this paper, it is important to discuss how the separation of substantive and statistical information can be achieved in order to make a case for treating statistical models as *canonical models* which can be used in conjunction with substantive information from any applied field.

It is widely recognized that stochastic phenomena amenable to empirical modeling have two *interrelated* sources of information, the *substantive subject matter* and the *statistical information* (chance regularity). What is not so apparent is how these sources of information are integrated in the context of empirical modeling. The PR approach treats the statistical and substantive information as complementary and, *ab initio,* are described separately in the form of a statistical and a structural model, respectively. The key for this *ab initio* separation is provided by viewing a statistical model generically as a particular parameterization of a stochastic processes $\{\mathbf{Z}_t, t \in \mathbb{T}\}$ underlying the data $\mathbf{Z}$, which, under certain conditions, can nest (parametrically) the structural model(s) in question. This gives rise to a framework for integrating the various facets of modeling encountered in the discussion of the early contributions by Fisher and Neyman: *specification, misspecification testing, respecification, statistical adequacy, statistical (inductive) inference, and identification.*

### 3.1. Structural vs. statistical models

It is widely recognized that most stochastic phenomena (the ones exhibiting chance regularity patterns) are commonly influenced by a very large number of contributing factors, and that explains why theories are often dominated by *ceteris paribus* clauses. The idea behind a *theory* is that in explaining the behavior of a variable, say $y_k$, one demarcates the segment of reality to be modeled by selecting the primary influencing factors $\mathbf{x}_k$, cognizant of the fact that there might be numerous other potentially relevant factors $\xi_k$ (observable and unobservable) that jointly determine the behavior of $y_k$ via a *theory model:*

$$(3.1) \qquad y_k = h^*(\mathbf{x}_k, \xi_k), \quad k \in \mathbb{N},$$

where $h^*(.)$ represents the *true* behavioral relationship for $y_k$. The guiding principle in selecting the variables in $\mathbf{x}_k$ is to ensure that they collectively account for the *systematic behavior* of $y_k$, and the *unaccounted factors* $\boldsymbol{\xi}_k$ represent non-essential disturbing influences which have only a *non-systematic* effect on $y_k$. This reasoning transforms (3.1) into a *structural model* of the form:

$$(3.2) \qquad y_k = h(\mathbf{x}_k; \boldsymbol{\phi}) + \epsilon(\mathbf{x}_k \boldsymbol{\xi}_k), \quad k \in \mathbb{N},$$

where $h(.)$ denotes the postulated functional form, $\boldsymbol{\phi}$ stands for the *structural parameters* of interest, and $\epsilon(\mathbf{x}_k \boldsymbol{\xi}_k)$ represents the *structural error term,* viewed as a function of both $\mathbf{x}_k$ and $\boldsymbol{\xi}_k$. By definition the error term process is:

$$(3.3) \qquad \{\epsilon(\mathbf{x}_k \boldsymbol{\xi}_k) = y_k - h(\mathbf{x}_k \boldsymbol{\phi}), k \in \mathbb{N}\},$$

and represents all unmodeled influences, *intended* to be a *white-noise* (non-systematic) process, i.e. for all possible values $(\mathbf{x}_k \boldsymbol{\xi}_k) \in \mathbb{R}_{\mathbf{x}} \times \mathbb{R}_{\boldsymbol{\xi}}$:

[i] $E[\epsilon(\mathbf{x}_k \boldsymbol{\xi}_k)] = 0$, [ii] $E[\epsilon(\mathbf{x}_k \boldsymbol{\xi}_k)^2] = \sigma_\epsilon^2$, [iii] $E[\epsilon(\mathbf{x}_k \boldsymbol{\xi}_k) \cdot \epsilon(\mathbf{x}_j \boldsymbol{\xi}_j)] = 0$, for $k \neq j$.

In addition, (3.2) represents a 'nearly isolated' generating mechanism in the sense that its error should be uncorrelated with the modeled influences (systematic component $h(\mathbf{x}_k \boldsymbol{\phi})$), i.e. [iv] $E[\epsilon(\mathbf{x}_k \boldsymbol{\xi}_k) \cdot h(\mathbf{x}_k \boldsymbol{\phi})] = 0$; the term 'nearly' refers to the non-deterministic nature of the isolation - see Spanos ([31], [35]).

In summary, a *structural model* provides an 'idealized' substantive description of the phenomenon of interest, in the form of a 'nearly isolated' mathematical system (3.2). The specification of a structural model comprises several choices: (a) the demarcation of the segment of the phenomenon of interest to be captured, (b) the important aspects of the phenomenon to be measured, and (c) the extent to which the inferences based on the structural model are germane to the phenomenon of interest. The kind of *errors* one can probe for in the context of a structural model concern the choices (a)–(c), including the form of $h(\mathbf{x}_k; \boldsymbol{\phi})$ and the circumstances that render the error term potentially systematic, such as the presence of relevant factors, say $\mathbf{w}_k$, in $\boldsymbol{\xi}_k$ that might have a systematic effect on the behavior of $y_t$; see Spanos [41].

It is important to emphasize that (3.2) depicts a 'factual' Generating Mechanism (GM), which aims to approximate the *actual data GM*. However, the assumptions [i]–[iv] of the structural error are *non-testable* because their assessment would involve verification for all possible values $(\mathbf{x}_k \boldsymbol{\xi}_k) \in \mathbb{R}_{\mathbf{x}} \times \mathbb{R}_{\boldsymbol{\xi}}$. To render them testable one needs to embed this structural into a statistical model; a crucial move that often goes unnoticed. Not surprisingly, the nature of the embedding itself depends crucially on whether the data $\mathbf{Z} := (\mathbf{z}_1, \mathbf{z}_2, \ldots, \mathbf{z}_n)$ are the result of an experiment or they are non-experimental (observational) in nature.

### 3.2. Statistical models and experimental data

In the case where one can perform experiments, 'experimental design' techniques, might allow one to operationalize the 'near isolation' condition (see Spanos, [35]), including the *ceteris paribus* clauses, and ensure that the *error term* is no longer a function of $(\mathbf{x}_k \boldsymbol{\xi}_k)$, but takes the generic form:

$$(3.4) \qquad \epsilon(\mathbf{x}_k \boldsymbol{\xi}_k) = \varepsilon_k \sim \mathsf{IID}(0, \sigma^2), \quad k = 1, 2, \ldots, n.$$

For instance, *randomization* and *blocking* are often used to 'neutralize' the phenomenon from the potential effects of $\boldsymbol{\xi}_k$ by ensuring that these uncontrolled factors cancel each other out; see Fisher [7]. As a direct result of the experimental 'control' via (3.4) the structural model (3.2) is essentially transformed into a *statistical model:*

$$(3.5) \qquad y_k = h(\mathbf{x}_k; \boldsymbol{\theta}) + \varepsilon_k, \quad \varepsilon_k \sim \mathsf{IID}(0, \sigma^2),\ k = 1, 2, \ldots, n.$$

The statistical error terms in (3.5) are qualitatively very different from the structural errors in (3.2) because they no longer depend on $(\mathbf{x}_k \boldsymbol{\xi}_k)$; the clause 'for all $(\mathbf{x}_k \boldsymbol{\xi}_k) \in \mathbb{R}_\mathbf{x} \times \mathbb{R}_{\boldsymbol{\xi}}$' has been rendered irrelevant. The most important aspect of embedding the structural (3.2) into the statistical model (3.5) is that, in contrast to [i]–[iv] for $\{\epsilon(\mathbf{x}_k \boldsymbol{\xi}_k),\ k \in \mathbb{N}\}$, the probabilistic assumptions $\mathsf{IID}(0, \sigma^2)$ concerning the *statistical error term* are rendered *testable*. That is, by operationalizing the 'near isolation' condition via (3.4), the error term has been tamed. For more precise inferences one needs to be more specific about the probabilistic assumptions defining the statistical model, including the functional form $h(.)$. This is because the more finical the probabilistic assumptions (the more constricting the statistical premises) the more precise the inferences; see Spanos [40].

The ontological status of the statistical model (3.5) is different from that of the structural model (3.2) in so far as (3.4) has operationalized the 'near isolation' condition. The statistical model has been 'created' as a result of the experimental design and control. As a consequence of (3.4) the *informational universe of discourse* for the statistical model (3.5) has been delimited to the probabilistic information relating to the observables $\mathbf{Z}_k$. This probabilistic structure, according to Kolmogorov's consistency theorem, can be fully described, under certain mild regularity conditions, in terms of the joint distribution $D(\mathbf{Z}_1, \mathbf{Z}_2, \ldots, \mathbf{Z}_n; \boldsymbol{\phi})$; see Doob [4]. It turns out that a statistical model can be viewed as a *parameterization* of the presumed probabilistic structure of the process $\{\mathbf{Z}_k,\ k \in \mathbb{N}\}$; see Spanos ([31], [35]).

In summary, a *statistical model* constitutes an 'idealized' *probabilistic description* of a stochastic process $\{\mathbf{Z}_k,\ k \in \mathbb{N}\}$, giving rise to data $\mathbf{Z}$, in the form of an internally consistent set of probabilistic assumptions, chosen to ensure that this data constitute a 'truly typical realization' of $\{\mathbf{Z}_k,\ k \in \mathbb{N}\}$.

In contrast to a *structural model*, once $\mathbf{Z}_k$ is chosen, a *statistical model* relies exclusively on the statistical information in $D(\mathbf{Z}_1, \mathbf{Z}_2, \ldots, \mathbf{Z}_n; \boldsymbol{\phi})$, that 'reflects' the chance regularity patterns exhibited by the data. Hence, a statistical model acquires 'a life of its own' in the sense that it constitutes a self-contained GM defined exclusively in terms of probabilistic assumptions pertaining to the observables $\mathbf{Z}_k := (y_k, \mathbf{X}_k)$. For example, in the case where $h(\mathbf{x}_k; \boldsymbol{\phi}) = \beta_0 + \boldsymbol{\beta}_1^\top \mathbf{x}_k$, and $\varepsilon_k \backsim \mathsf{N}(.,.)$, (3.5) becomes the *Gauss Linear model,* comprising the statistical GM:

$$(3.6) \qquad y_k = \beta_0 + \boldsymbol{\beta}_1^\top \mathbf{x}_k + u_k,\ k \in \mathbb{N},$$

together with the probabilistic assumptions (Spanos [31]):

$$(3.7) \qquad y_k \backsim \mathsf{NI}(\beta_0 + \boldsymbol{\beta}_1^\top \mathbf{x}_k,\ \sigma^2),\ k \in \mathbb{N},$$

where $\boldsymbol{\theta} := (\beta_0, \boldsymbol{\beta}_1, \sigma^2)$ is assumed to be $k$-invariant, and 'NI' stands for 'Normal, Independent'.

### 3.3. Statistical models and observational data

This is the case where the observed data on $(y_t, \mathbf{x}_t)$ are the result of an ongoing actual data generating process, undisturbed by any experimental control or intervention. In this case the route followed in (3.4) in order to render the *statistical error term* (a) free of $(\mathbf{x}_t, \xi_t)$, and (b) non-systematic in a statistical sense, is no longer feasible. It turns out that sequential *conditioning* supplies the primary tool in modeling observational data because it provides an alternative way to ensure the non-systematic nature of the statistical error term without controls and intervention.

It is well-known that sequential conditioning provides a general way to transform an arbitrary stochastic process $\{\mathbf{Z}_t,\ t \in \mathbb{T}\}$ into a *Martingale Difference (MD) process* relative to an increasing sequence of sigma-fields $\{\mathfrak{D}_t,\ t \in \mathbb{T}\}$; a modern form of a non-systematic error process (Doob, [4]). This provides the key to an alternative approach to specifying statistical models in the case of non-experimental data by replacing the 'controls' and 'interventions' with the choice of the *relevant conditioning information set* $\mathfrak{D}_t$ that would render the error term a MD; see Spanos [31].

As in the case of experimental data the universe of discourse for a statistical model is fully described by the joint distribution $D(\mathbf{Z}_1, \mathbf{Z}_2, \ldots, \mathbf{Z}_T; \boldsymbol{\phi})$, $\mathbf{Z}_t := (y_t, \mathbf{X}_t^\top)^\top$. Assuming that $\{\mathbf{Z}_t,\ t \in \mathbb{T}\}$ has bounded moments up to order two, one can choose the conditioning information set to be:

(3.8) $$\mathfrak{D}_{t-1} = \sigma\left(y_{t-1}, y_{t-2}, \ldots, y_1, \mathbf{X}_t, \mathbf{X}_{t-1}, \ldots, \mathbf{X}_1\right).$$

This renders the error process $\{u_t,\ t \in \mathbb{T}\}$, defined by:

(3.9) $$u_t = y_t - E(y_t | \mathfrak{D}_{t-1}),$$

a MD process relative to $\mathfrak{D}_{t-1}$, irrespective of the probabilistic structure of $\{\mathbf{Z}_t,\ t \in \mathbb{T}\}$; see Spanos [36]. This error process is based on $D(y_t | \mathbf{X}_t, \mathbf{Z}_{t-1}^0; \psi_{1t})$, where $\mathbf{Z}_{t-1}^0 := (\mathbf{Z}_{t-1}, \ldots, \mathbf{Z}_1)$, which is directly related to $D(\mathbf{Z}_1, \ldots, \mathbf{Z}_T; \boldsymbol{\phi})$ via:

(3.10) $$\begin{aligned} D(\mathbf{Z}_1, \ldots, \mathbf{Z}_T; \boldsymbol{\phi}) \\ = D(\mathbf{Z}_1; \psi_1) \prod_{t=2}^{T} D_t(\mathbf{Z}_t | \mathbf{Z}_{t-1}^0; \psi_t) \\ = D(\mathbf{Z}_1; \psi_1) \prod_{t=2}^{T} D_t(y_t | \mathbf{X}_t, \mathbf{Z}_{t-1}^0; \psi_{1t}) \cdot D_t(\mathbf{X}_t | \mathbf{Z}_{t-1}^0; \psi_{2t}). \end{aligned}$$

The Greek letters $\phi$ and $\psi$ are used to denote the unknown parameters of the distribution in question. This sequential conditioning gives rise to a statistical GM of the form:

(3.11) $$y_t = E(y_t | \mathfrak{D}_{t-1}) + u_t,\ t \in \mathbb{T},$$

which is *non-operational* as it stands because without further restrictions on the process $\{\mathbf{Z}_t,\ t \in \mathbb{T}\}$, the systematic component $E(y_t | \mathfrak{D}_{t-1})$ cannot be specified explicitly. For operational models one needs to postulate some probabilistic structure for $\{\mathbf{Z}_t,\ t \in \mathbb{T}\}$ that would render the data $\mathbf{Z}$ a 'truly typical' realization thereof. These assumptions come from a menu of three broad categories: (D) **Distribution**, (M) **Dependence**, (H) **Heterogeneity**; see Spanos ([34]–[38]).

**Example.** *The Normal/Linear Regression model* results from the reduction (3.10) by assuming that $\{\mathbf{Z}_t,\ t \in \mathbb{T}\}$ is a NIID vector process. These assumptions ensure that the relevant information set that would render the error process a MD is reduced from $\mathfrak{D}_{t-1}$ to $\mathfrak{D}_t^x=\{\mathbf{X}_t=\mathbf{x}_t\}$, ensuring that:

$$(3.12) \qquad (u_t|\,\mathbf{X}_t=\mathbf{x}_t) \sim \mathsf{NIID}(0,\sigma^2),\ k=1,2,\ldots,T.$$

This is analogous to (3.4) in the case of experimental data, but now the error term has been operationalized by a judicious choice of $\mathfrak{D}_t^x$. The Linear Regression model comprises the statistical GM:

$$(3.13) \qquad y_t = \beta_0 + \boldsymbol{\beta}_1^\top \mathbf{x}_t + u_t,\ t \in \mathbb{T},$$

$$(3.14) \qquad (y_t|\,\mathbf{X}_t=\mathbf{x}_t) \sim \mathsf{NI}(\beta_0 + \boldsymbol{\beta}_1^\top \mathbf{x}_t,\ \sigma^2),\ t \in \mathbb{T},$$

where $\boldsymbol{\theta} := (\beta_0, \boldsymbol{\beta}_1, \sigma^2)$ is assumed to be $t$-invariant; see Spanos [35].

The probabilistic perspective gives a statistical model 'a life of its own' in the sense that the probabilistic assumptions in (3.14) bring to the table *statistical information* which supplements, and can be used to assess the appropriateness of, the *substantive subject matter information*. For instance, in the context of the structural model $h(\mathbf{x}_t; \boldsymbol{\phi})$ is determined by the theory. In contrast, in the context of a statistical model it is determined by the probabilistic structure of the process $\{\mathbf{Z}_t,\ t\in\mathbb{T}\}$ via $h(\mathbf{x}_t; \boldsymbol{\theta})=E(y_t|\,\mathbf{X}_t=\mathbf{x}_t)$, which, in turn, is determined by the joint distribution $D(y_t, \mathbf{X}_t; \psi)$; see Spanos [36].

An important aspect of embedding a structural into a statistical model is to ensure (whenever possible) that the former can be viewed as a *reparameterization/restriction* of the latter. The structural model is then tested against the benchmark provided by a statistically adequate model. *Identification* refers to being able to define $\boldsymbol{\phi}$ uniquely in terms of $\boldsymbol{\theta}$. Often $\boldsymbol{\theta}$ has more parameters than $\boldsymbol{\phi}$ and the embedding enables one to test the validity of the additional restrictions, known as over-identifying restrictions; see Spanos [33].

### 3.4. Kepler's first law of planetary motion revisited

In an attempt to illustrate some of the concepts and procedures introduced in the PR framework, we revisit Lehmann's [11] example of Kepler's statistical model predating, by more than 60 years, the eventual structural model proposed by Newton.

Kepler's law of planetary motion was originally just an *empirical regularity* that he 'deduced' from Brahe's data, stating that the motion of any planet around the sun is elliptical. That is, the loci of the motion in *polar coordinates* takes the form $(1/r)=\alpha_0+\alpha_1\cos\vartheta$, where $r$ denotes *the distance of the planet from the sun*, and $\vartheta$ denotes *the angle between the line joining the sun and the planet and the principal axis of the ellipse*. Defining the observable variables by $y:=(1/r)$ and $x:=\cos\vartheta$, Kepler's empirical regularity amounted to an estimated linear regression model:

$$(3.15) \qquad y_t = \underset{(.000002)}{0.662062} + \underset{(.000003)}{.061333} x_t + \widehat{u}_t,\ R^2=.999,\ s=.0000111479;$$

these estimates are based on Kepler's original 1609 data on Mars with $n=28$. Formal misspecification tests of the model assumptions in (3.14) (Section 3.3), indicate that the estimated model is statistically adequate; see Spanos [39] for the details.

Substantive interpretation was bestowed on (3.15) by *Newton's law of universal gravitation*: $F=\frac{G(m \cdot M)}{r^2}$, where $F$ is the force of attraction between two bodies of *mass m* (planet) and $M$ (sun), $G$ is a *constant of gravitational attraction*, and $r$ is the *distance* between the two bodies, in the form of a *structural model*:

(3.16) $$Y_k = \alpha_0 + \alpha_1 X_k + \epsilon(x_k, \xi_k), \ k \in \mathbb{N},$$

where the parameters $(\alpha_0, \alpha_1)$ are given a *structural interpretation*: $\alpha_0 = \frac{MG}{4\kappa^2}$, where $\kappa$ denotes Kepler's constant, $\alpha_1 = (\frac{1}{d} - \alpha_0)$, $d$ denotes the shortest distance between the planet and the sun. The error term $\epsilon(x_k, \xi_k)$ also enjoys a structural interpretation in the form of *unmodeled effects*; its assumptions [i]–[iv] (Section 3.1) will be *inappropriate* in cases where (a) the data suffer from 'systematic' observation errors, and there are significant (b) *third body* and/or (c) *general relativity* effects.

### 3.5. Revisiting certain issues in empirical modeling

In what follows we indicate very briefly how the PR approach can be used to shed light on certain crucial issues raised by Lehmann [11] and Cox [1].

**Specification: a 'Fountain' of statistical models.** The PR approach broadens Lehmann's *reservoir of models* idea to the set of all possible statistical models $\mathcal{P}$ that could (potentially) have given rise to data **Z**. The statistical models in $\mathcal{P}$ are *characterized* by their reduction assumptions from three broad categories: *Distribution, Dependence,* and *Heterogeneity*. This way of viewing statistical models provides (i) a systematic way to *characterize statistical models*, (different from Lehmann's) and, at the same time it offers (ii) a general procedure to generate new statistical models.

The capacity of the PR approach to generate new statistical models is demonstrated in Spanos [36], ch. 7, were several bivariate distributions are used to derive different regression models via (3.10); this gives rise to several non-linear and/or heteroskedastic regression models, most of which remain unexplored. In the same vein, the reduction assumptions of (D) *Normality,* (M) *Markov dependence,* and (H) *Stationarity*, give rise to Autoregressive models; see Spanos ([36], [38]).

Spanos [34] derives a new family of Linear/heteroskedastic regression models by replacing the Normal in (3.10) with the Student's t distribution. When the IID assumptions were also replaced by Markov dependence and Stationarity, a surprising family of models emerges that extends the ARCH formulation; see McGuirk et al [15], Heracleous and Spanos [8].

**Model validation: statistical vs. structural adequacy.** The PR approach also addresses Lehmann's concern that structural and statistical models 'pose very different problems for model validation'; see Spanos [41]. The purely probabilistic construal of statistical models renders *statistical adequacy* the only relevant criterion for model validity is *statistical adequacy*. This is achieved by thorough *misspecification testing* and *respecification*; see Mayo and Spanos [13].

MisSpecification (M-S) testing is different from Neyman and Pearson (N–P) testing in one important respect. N–P testing assumes that the prespecified statistical model class $\mathcal{M}$ includes the true model, say $f_0(\mathbf{z})$, and probes *within the boundaries* of this model using the hypotheses:

$$H_0: f_0(\mathbf{z}) \in \mathcal{M}_0 \text{ vs. } H_1: f_0(\mathbf{z}) \in \mathcal{M}_1,$$

where $\mathcal{M}_0$ and $\mathcal{M}_1$ form a partition of $\mathcal{M}$. In contrast, M-S testing probes *outside the boundaries* of the prespecified model:

$$H_0: f_0(\mathbf{z}) \in \mathcal{M} \text{ vs. } \overline{H}_0: f_0(\mathbf{z}) \in [\mathcal{P} - \mathcal{M}],$$

where $\mathcal{P}$ denotes the set of all possible statistical models, rendering them Fisherian type *significance tests.* The problem is how one can operationalize $\mathcal{P} - \mathcal{M}$ in order to probe thoroughly for possible departures; see Spanos [36]. Detection of departures from the null in the direction of, say $\mathcal{P}_1 \subset [\mathcal{P} - \mathcal{M}]$, is sufficient to deduce that the null is false but not to deduce that $\mathcal{P}_1$ is true; see Spanos [37]. More formally, $\mathcal{P}_1$ *has not passed a severe test*, since its own statistical adequacy has not been established; see Mayo and Spanos ([13], [14]).

On the other hand, validity for a structural model refers to *substantive adequacy*: a combination of *data-acceptability* on the basis of a statistically adequate model, and *external validity* - how well the structural model 'approximates' the reality it aims to explain. Statistical adequacy is a precondition for the assessment of substantive adequacy because without it no reliable inference procedures can be used to assess substantive adequacy; see Spanos [41].

**Model specification vs. model selection**. The PR approach can shed light on Lehmann's concern about *model specification* vs. *model selection,* by underscoring the fact that the primary criterion for model specification within $\mathcal{P}$ is statistical adequacy, not goodness of fit. As pointed out by Lehmann [11], the current *model selection* procedures (see Rao and Wu, [30], for a recent survey) do not address the original statistical model specification problem. One can make a strong case that Akaike-type model selection procedures assume the statistical model specification problem solved. Moreover, when the statistical adequacy issue is addressed, these model selection procedure becomes superfluous; see Spanos [42].

**Statistical Generating Mechanism (GM)**. It is well-known that a statistical model can be specified fully in terms of the joint distribution of the observable random variables involved. However, if the statistical model is to be related to any structural models, it is imperative to be able to specify a statistical GM which will provide the bridge between the two models. This is succinctly articulated by Cox [1]:

"The essential idea is that if the investigator cannot use the model directly to simulate artificial data, how can "Nature" have used anything like that method to generate real data?" (p. 172)

The PR specification of statistical models brings the statistical GM based on the orthogonal decomposition $y_t = E(y_t | \mathfrak{D}_{t-1}) + u_t$ in (3.11) to the forefront. The onus is on the modeler to choose (i) an appropriate probabilistic structure for $\{y_t,\ t \in \mathbb{T}\}$, and (ii) the associated information set $\mathfrak{D}_{t-1}$, relative to which the error term is rendered a martingale difference (MD) process; see Spanos [36].

**The role of exploratory data analysis**. An important feature of the PR approach is to render the use of graphical techniques and exploratory data analysis (EDA), more generally, an integral part of statistical modeling. EDA plays a crucial role in the specification, M-S testing and respecification facets of modeling. This addresses a concern raised by Cox [1] that:

"... the separation of 'exploratory data analysis' from 'statistics' are counterproductive." (ibid., p. 169)

## 4. Conclusion

Lehmann [11] raised the question whether the presence of substantive information subordinates statistical modeling to other disciplines, precluding statistics from having its own intended scope. This paper argues that, despite the uniqueness of every modeling endeavor arising from the *substantive subject matter information*, all forms of statistical modeling share certain generic aspects which revolve around the notion of *statistical information*. The key to upholding the integrity of both sources of information, as well as ensuring the reliability of their fusing, is a purely probabilistic construal of statistical models in the spirit of Fisher and Neyman. The PR approach adopts this view of specification and accommodates the related facets of modeling: *misspecification testing* and *respecification*.

The PR modeling framework gives the statistician a *pivotal role* and extends the intended scope of statistics, without relegating the role of substantive information in empiridal modeling. The judicious use of probability theory, in conjunction with graphical techniques, can transform the specification of statistical models into purpose-built conjecturing which can be assessed subsequently. In addition, thorough misspecification testing can be used to assess the appropriateness of a statistical model, in order to ensure the reliability of inductive inferences based upon it. Statistically adequate models *have a life of their own* in so far as they can be (sometimes) the ultimate objective of modeling or they can be used to establish *empirical regularities* for which substantive explanations need to account; see Cox [1]. By embedding a structural into a statistically adequate model and securing substantive adequacy, confers upon the former statistical meaning and upon the latter substantive meaning, rendering learning from data, using statistical induction, a reliable process.